\documentclass[11pt, a4paper]{amsart}

\usepackage[all]{xy}
\usepackage{amsmath}
\usepackage{amssymb}
\usepackage{amsfonts}
\usepackage{graphicx}
\usepackage[mathscr]{eucal}
\usepackage{enumerate}
\usepackage[latin1]{inputenc}
\usepackage{times}

\theoremstyle{definition}
\newtheorem{definition}{Definition}[section]
\theoremstyle{remark}

\newtheoremstyle{estilo}{\topsep}{\topsep}{\slshape}{}{\bfseries}{.}{ }{}
\theoremstyle{estilo}
\newtheorem{theorem}[definition]{Theorem}
\newtheorem{corollary}[definition]{Corollary}
\newtheorem{lemma}[definition]{Lemma}
\newtheorem{proposition}[definition]{Proposition}
\newtheorem{remark}[definition]{Remark}

\newtheorem{examples}[definition]{Examples}

\def\ot{\otimes}

\newcommand{\nat}{\mbox{$\;\natural \;$}}

\renewenvironment{proof}{{\noindent\sc Proof\;}}{\qed\\}

\def\O{\Omega}
\def\c{\mathbb C}
\newcommand{\dtext}[1]{\emph{\textbf{#1}}}
\def\vphi{\varphi}
\def\eps{\varepsilon}
\newcommand{\abs}[1]{\left\lvert#1\right\rvert}
\def\M{{\mathcal{M}}}
\def\R{\mathbb R}

\def\l{\lambda}

\newcommand{\LTwistor}[1]
    {\xy 
    		(4,6)*{\boxed{\scriptscriptstyle{#1}}}="x";
      	(0,12)*{}; "x" **\crv{(0,8)&(2,12)};
      	"x"; (8,0)*{} **\crv{(6,0)&(8,4)}; 
        (0,0)*{}; "x" **\crv{(0,4)&(2,0)};
        (8,12)*{}; "x" **\crv{(8,8)&(6,12)}; 
    \endxy}

\begin{document}

\title{General twisting of algebras}
\thanks{Research partially supported by the EC programme LIEGRITS,
RTN 2003, 505078. The first author has been also partially supported
by the projects MTM2004-08125 and FQM-266 (Junta de Andaluc\'{i}a
Research Group) and Spanish MEC-FPU grant AP2003-4340. The second
and the third authors have been also partially supported by the
bilateral project ``New techniques in Hopf algebras and graded ring
theory'' of the Flemish and Romanian Ministries of Research. The
second author has been also partially supported by the CEEX
programme of the Romanian Ministry of Education and Research,
contract nr. CEx05-D11-11/2005.}
\author[J. L\'{o}pez]{Javier L\'{o}pez Pe\~{n}a}
\address{Department of Algebra,
University of Granada\\
Avda. Fuentenueva s/n, E-18071, Granada, Spain}
\email{jlopez@ugr.es}
\author[F. Panaite]{Florin Panaite}
\address{Institute of Mathematics,
Romanian Academy\\
PO-Box 1-764, RO-014700 Bucharest, Romania}
\email{Florin.Panaite@imar.ro}
\author[F. Van Oystaeyen]{Freddy Van
Oystaeyen}
\address{Department of Mathematics and Computer Science\\
University of Antwerp, Middelheimlaan 1\\
B-2020 Antwerp, Belgium}
\email{Francine.Schoeters@ua.ac.be
(secretary)}
\date{}
\begin{abstract}
We introduce the concept of {\it pseudotwistor} (with particular cases
called {\it twistor} and {\it braided twistor}) for an algebra
$(A, \mu , u)$ in a monoidal category, as a morphism
$T:A\otimes A\rightarrow A\otimes A$ satisfying a list of axioms
ensuring that $(A, \mu \circ T, u)$ is also an algebra in the category.
This concept provides a unifying framework for various deformed (or twisted)
algebras from the literature, such as twisted tensor products of algebras,
twisted bialgebras and algebras endowed with Fedosov products. Pseudotwistors
appear also in other topics from the literature, e.g. 
Durdevich's braided quantum
groups and ribbon algebras.
We also focus on the effect of twistors on the universal first 
order differential
calculus, as well as on lifting twistors to braided twistors on 
the algebra of universal
differential forms.
\end{abstract}
\maketitle
\section{Introduction}
The twisted tensor product $A\otimes _R B$ of two associative
algebras $A$ and $B$ is a certain associative algebra structure on
the vector space $A\otimes B$, defined in terms of a so-called
{\it twisting map} $R:B\otimes A \rightarrow A\otimes B$, having the 
property that it coincides with the usual tensor product $A\otimes
B$ if $R$ is the usual flip. This construction was proposed in
\cite{cap} as a representative for the cartesian product of
noncommutative spaces. More evidence that this proposal is
meaningful appeared recently in \cite{jlpvo}, where it was proved
that this construction may be iterated in a natural way, and that
the noncommutative $2n$-planes defined by Connes and
Dubois-Violette, cf. \cite{cdv}, may be written as iterated twisted
tensor products of some commutative algebras. Various other
applications of twisted tensor products appear in the literature,
see for instance \cite{bor}, \cite{vandaele}. Note also that, as we 
learned from the referee, categorical analogues of twisting maps 
appeared earlier in the literature, under the name {\it distributive laws}, 
see for instance \cite{beck}, \cite{markl}, \cite{street}.

On the other hand, if $H$ is a bialgebra and $\sigma :H\otimes
H\rightarrow k$ is a normalized and convolution invertible left
2-cocycle, one can consider the ``twisted bialgebra'' $_{\sigma }H$,
which is an associative algebra structure on $H$ with multiplication
$a*b=\sigma (a_1, b_1)a_2b_2$. This is an important and well-known
construction, containing as particular case the classical twisted
group rings.

Apparently, there is no relation between twisted tensor products of
algebras and twisted bialgebras, except for the fact that their
names suggest that they are both obtained via a process of {\it
twisting}. However, as a consequence of the ideas developed in this
paper, it will turn out that this suggestion is correct: we will
find a framework in which both these constructions fit as particular
cases.

 Our initial aim was to relate the multiplications $\mu
_{A\otimes _R B}$ of $A\otimes _R B$ and $\mu _{A\otimes B}$ of
$A\otimes B$. It is easy to see that $\mu _{A\otimes _R B}=\mu
_{A\otimes B}\circ T$, where $T:(A\otimes B)\otimes (A\otimes
B)\rightarrow (A\otimes B)\otimes (A\otimes B)$ is a map depending
on $R$, and the problem is to find the abstract properties satisfied
by this map $T$, which together with the associativity of $\mu
_{A\otimes B}$ imply the associativity of $\mu _{A\otimes _R B}$. We
are thus led to introduce the concept of {\it twistor} for an
algebra $D$, as a linear map $T:D\otimes D\rightarrow D\otimes D$
satisfying a list of axioms which imply that the new multiplication
$\mu _D\circ T$ is an associative algebra structure on the vector
space $D$ (these axioms are similar to, but different from, the ones
of an $R$-matrix for an associative algebra, a concept introduced by
Borcherds). It turns out that the map $T$ affording the
multiplication of $A\otimes _R B$ is such a twistor, and that
various other examples of twistors may be identified in the
literature, in particular the noncommutative $2n$-plane may be
regarded as a deformation of a polynomial algebra via a twistor.

But there exist in the literature many examples of deformed
multiplications which are {\it not} afforded by twistors. For
instance, the map $T(a\otimes b)=\sigma (a_1, b_1)a_2\otimes b_2$
affording the multiplication of $_{\sigma }H$ is far from being a
twistor. But the map $T(\omega \otimes \zeta )=\omega \otimes \zeta
-(-1)^{|\omega |} d(\omega )\otimes d(\zeta )$, affording the
so-called Fedosov product, is not too far, it looks like a graded
analogue. We are thus led to a more general concept, called {\it
braided twistor}, of which this $T$ is an example. And from this
concept we arrive at a much more general one, called {\it
pseudotwistor}, which is general enough to include as example the
map affording the multiplication of $_{\sigma }H$, as well as some
other (nonrelated) situations from the literature, e.g. some
examples arising in the context of Durdevich's braided quantum
groups, and the morphism $c_{A, A}^2$, where $A$ is an algebra in a
braided monoidal category with braiding $c$.

 We also present some
properties of (pseudo)twistors, e.g. we show how to lift modules and
bimodules over $D$ to the same structures over the deformed algebra,
and how to extend a twistor $T$ for an algebra $D$ to a braided
(graded) twistor $\widetilde{T}$ for the algebra of universal
differential forms $\O D$.
\section{Preliminaries}
\setcounter{equation}{0}
Let $k$ be a field, used as a base field throughout. We denote
$\otimes _k$ by $\otimes $, the identity $id_V$ of an object $V$
simply by $V$, and by $\tau :V\otimes W \rightarrow W\otimes V$,
$\tau (v\otimes w)=w\otimes v$, the usual flip. All algebras are
assumed to be associative unital $k$-algebras; the multiplication
and unit of an algebra $D$ are denoted by $\mu _D:D\otimes
D\rightarrow D$ and respectively $u_D:k\rightarrow D$ (or simply by
$\mu $ and $u$ if there is no danger of confusion). For bialgebras
and Hopf algebras we use the Sweedler-type notation $\Delta
(h)=h_1\otimes h_2$, and for categorical terminology we refer to
\cite{joyal}, \cite{k}, \cite{m2}. For some proofs, we will use
braiding notation, of which a detailed description may be found in
\cite{k}.

 We recall the twisted tensor product of algebras from
\cite{tambara}, \cite{vandaele}, \cite{cap}. If $A$ and $B$ are two
algebras, a linear map $R:B\ot A\rightarrow A\ot B$ is called a
\dtext{twisting map} if it satisfies the conditions
\begin{eqnarray}
&&R(b\ot 1)=1\ot b,\;\;\;R(1\ot a)=a\ot 1, \;\;\;\forall \; a\in A,
\;b\in B, \label{tw1}\\
&&R\circ (B\ot \mu _A)=(\mu _A\ot B)\circ (A\ot R)\circ (R\ot A),
\label{tw2}\\
&&R\circ (\mu _B\ot A)=(A\ot \mu _B)\circ (R\ot B)\circ (B\ot R).
\label{tw3}
\end{eqnarray}
If we denote by $R(b\ot a)=a_R\ot b_R$, for $a\in A$, $b\in B$, then
(\ref{tw2}) and (\ref{tw3}) may be written as:
\begin{eqnarray}
&&(aa')_R\ot b_R=a_Ra'_r\ot (b_R)_r, \label{tw4} \\
&&a_R\ot (bb')_R=(a_R)_r\ot b_rb'_R, \label{tw5}
\end{eqnarray}
for all $a, a'\in A$ and $b, b'\in B$, where $r$ is another copy of $R$.
If we define a multiplication on $A\ot B$, by $\mu _R=(\mu _A\otimes
\mu _B)\circ (A\otimes R\otimes B)$, that is
\begin{eqnarray}
&&(a\ot b)(a'\ot b')=aa'_R\ot b_Rb', \label{multtw}
\end{eqnarray}
then this multiplication is associative and $1\ot 1$ is the unit.
This algebra structure is denoted by $A\ot _RB$ and is called the
\dtext{twisted tensor product} of $A$ and $B$. This construction
works also if $A$ and $B$ are algebras in an arbitrary monoidal
category.

 If $A\otimes _{R_1} B$, $B\otimes _{R_2} C$ and $A\otimes
_{R_3} C$ are twisted tensor products of algebras, the twisting maps
$R_1$, $R_2$, $R_3$ are called {\it compatible} if they satisfy
\begin{eqnarray*}
&&(A\otimes R_2)\circ (R_3\otimes B)\circ (C\otimes R_1)=
(R_1\otimes C)\circ (B\otimes R_3)\circ (R_2\otimes A),
\end{eqnarray*}
see \cite{jlpvo}. If this is the case, the maps $T_1:C\otimes
(A\otimes _{R_1} B)\rightarrow (A\otimes _{R_1} B)\otimes C$ and
$T_2:(B\otimes _{R_2} C)\otimes A\rightarrow A\otimes (B\otimes
_{R_2} C)$ given by $T_1:=(A\otimes R_2)\circ (R_3\otimes B)$ and
$T_2:=(R_1\otimes C)\circ (B\otimes R_3)$ are also twisting maps and
$A\otimes _{T_2} (B\otimes _{R_2} C)\equiv (A\otimes _{R_1} B)
\otimes _{T_1} C$; this algebra is denoted by $A\otimes _{R_1} B
\otimes _{R_2} C$. This construction may be iterated to an arbitrary
number of factors, see \cite{jlpvo} for complete detail.

We recall the following result from \cite{cap}, to be used in the
sequel:
\begin{theorem}\label{twistdifferentialforms}
Let $A$, $B$ be two algebras. Then any twisting map $R:B\otimes A\to
A\otimes B$ extends to a unique twisting map $\tilde{R}:\O B\otimes
\O A \to \O A\otimes \O B$ which satisfies the conditions
\begin{eqnarray}\label{twisteddiff1}
\tilde{R}\circ(d_B\otimes \O A)&=&(\eps_A \otimes d_B)\circ
\tilde{R}, \\
\label{twisteddiff2} \tilde{R}\circ (\O B\otimes d_A)&=&(d_A \otimes
\eps_B)\circ \tilde{R},
\end{eqnarray}
where $d_A$ and $d_B$ denote the differentials on the algebras of 
universal differential forms $\O A$ and $\O B$, 
and $\eps_A$, $\eps_B$ stand for the gradings on $\O A$ and $\O B$,
respectively. Moreover, $\O A\otimes_{\tilde{R}}\O B$ is a graded
differential algebra with differential $d(\vphi\otimes
\omega):=d_A\vphi\otimes \omega + (-1)^{\abs{\vphi}}\vphi\otimes
d_B\omega$.
\end{theorem}

Finally, we recall the definition of the noncommutative $2n$-planes
introduced by Connes and Dubois-Violette in \cite{cdv}. Consider
$\theta\in \M_n(\R)$ an antisymmetric matrix,
$\theta=(\theta_{\mu\nu})$, $\theta_{\nu\mu}=-\theta_{\mu\nu}$, and
let $C_{alg}(\R^{2n}_\theta)$ be the associative algebra generated
by $2n$ elements $\{z^\mu,\bar{z}^\mu\}_{\mu=1,\dotsc,n}$ with
relations
\begin{equation} \label{ncplanerelations}
\left.\begin{array}{r}
z^\mu z^\nu=\l^{\mu\nu}z^\nu z^\mu \\
\bar{z}^\mu \bar{z}^\nu=\l^{\mu\nu}\bar{z}^\nu \bar{z}^\mu \\
\bar{z}^\mu z^\nu = \l^{\nu \mu} z^\nu \bar{z}^\mu
\end{array}\right\}\forall\,\mu,\nu=1,\dotsc,n,\ \text{being
$\l^{\mu\nu}:=e^{i\theta_{\mu\nu}}$.}
\end{equation}
Note that $\l^{\nu\mu}=(\l^{\mu\nu})^{-1}=\overline{\l^{\mu\nu}}$
for $\mu\neq \nu$, and $\l^{\mu\mu}=1$ by antisymmetry.
The algebra
$C_{alg}(\R^{2n}_\theta)$ will be then referred to as the
\dtext{(algebra of complex polynomial functions on the)
noncommutative $2n$--plane $\R^{2n}_\theta$}. In fact, former
relations define a deformation $\c^n_\theta$ of $\c^n$, so we can
identify the noncommutative complex $n$--plane $\c^n_\theta$ with
$\R^{2n}_\theta$ by writing
$C_{alg}(\c^n_\theta):=C_{alg}(\R^{2n}_\theta)$. As shown in  \cite{jlpvo},
$C_{alg}(\R^{2n}_\theta)$ may be written as an iterated twisted tensor product
of $n$ commutative (polynomial) algebras.

\section{$R$-matrices and twistors}\label{rmattw}
\setcounter{equation}{0}
In the literature there exist various schemes producing, from
a given associative algebra $A$ and some datum corresponding to it,
a new associative algebra structure on the vector space $A$.
The aim of this section is to prove that there exists such a
general scheme that produces the twisted tensor product starting from
the ordinary tensor product. Our source of inspiration is the following
result of Borcherds from \cite{borcherds1}, \cite{borcherds2},
which arose in his Hopf algebraic approach to vertex algebras:
\begin{theorem} (\cite{borcherds1}, \cite{borcherds2}) \label{borch}
Let $D$ be an algebra with multiplication denoted by
$\mu _D=\mu $ and let $T:D\ot D\rightarrow D\ot D$ be a linear map
satisfying the
following conditions: $T(1\ot d)=1\ot d$, $T(d\ot 1)=d\ot 1$, for all
$d\in D$, and
\begin{eqnarray}
&&\mu _{23}\circ T_{12}\circ T_{13}=
T\circ \mu _{23}:D\ot D\ot D\rightarrow D\ot D, \label{rmat1}\\
&&\mu _{12}\circ T_{23}\circ T_{13}=
T\circ \mu _{12}:D\ot D\ot D\rightarrow D\ot D, \label{rmat2}\\
&&T_{12}\circ T_{13}\circ T_{23}=
T_{23}\circ T_{13}\circ T_{12}:D\ot D\ot D\rightarrow D\ot D\ot D,
\label{rmat3}
\end{eqnarray}
with standard notation for $\mu _{ij}$ and $T_{ij}$. Then the
bilinear map $\mu \circ T:D\ot D\rightarrow D$ is another
associative algebra structure on $D$, with the same unit 1. The map
$T$ is called an \dtext{$R$-matrix}.
\end{theorem}

If $A\ot _RB$ is a twisted tensor product of algebras, we want to obtain it
as a twisting of $A\ot B$. Define
$T:(A\ot B)\ot (A\ot B)\rightarrow (A\ot B)\ot (A\ot B)$ by
$T=(A\otimes \tau \otimes B)\circ (A\otimes R\otimes B)$, i. e.
\begin{eqnarray}
T((a\ot b)\ot (a'\ot b'))=(a\ot b_R)\ot (a'_R\ot b'). \label{t}
\end{eqnarray}
Then the multiplication of $A\ot _RB$ is obtained as $\mu _{A\ot B}\circ T$,
also
$T$ satisfies $T(1\ot (a\ot b))=1\ot (a\ot b)$ and $T((a\ot b)\ot 1)=
(a\ot b)\ot 1$, but in general $T$ does {\it not} satisfy the other axioms in
Theorem \ref{borch} (for instance take $R$ to be the twisting map
corresponding to a Hopf smash product),
hence we {\it cannot} obtain $A\ot _RB$ from
$A\ot B$ using Borcherds' scheme, we have to find an alternative one.
This is achieved in the next result (the proof is postponed to Section
\ref{secbraided}, where it will be given in a more general framework).

\begin{theorem} \label{scheme}
Let $D$ be an algebra with multiplication denoted by
$\mu _D=\mu $ and $T:D\ot D\rightarrow D\ot D$ a linear map satisfying the
following conditions: $T(1\ot d)=1\ot d$, $T(d\ot 1)=d\ot 1$, for all
$d\in D$, and
\begin{eqnarray}
&&\mu _{23}\circ T_{13}\circ T_{12}=
T\circ \mu _{23}:D\ot D\ot D\rightarrow D\ot D, \label{dec1} \\
&&\mu _{12}\circ T_{13}\circ T_{23}=
T\circ \mu _{12}:D\ot D\ot D\rightarrow D\ot D,  \label{dec2}\\
&&T_{12}\circ T_{23}=T_{23}\circ T_{12}:
D\ot D\ot D\rightarrow D\ot D\ot D. \label{dec3}
\end{eqnarray}
Then the bilinear map $\mu \circ T:D\ot D\rightarrow D$ is another
associative algebra structure on $D$, with the same unit 1, which
will be denoted in what follows by $D^T$, and the map $T$ will be
called a \dtext{twistor} for $D$.
\end{theorem}

If $T$ is a twistor, we will usually denote $T(d\ot d')=d^T\ot d'_T$,
for $d, d'\in D$, so the new multiplication $\mu \circ T$
on $D$ is given by $d*d'=d^Td'_T$. With this notation, the relations
(\ref{dec1})--(\ref{dec3}) may be written as:
\begin{eqnarray}
&&d^T\ot (d'd'')_T=(d^T)^t\ot d'_Td''_t, \label{dec4} \\
&&(dd')^T\ot d''_T=d^Td'^t\ot (d''_t)_T, \label{dec5} \\
&&d^T\ot (d'_T)^t\ot d''_t=d^T\ot (d'^t)_T\ot d''_t. \label{dec6}
\end{eqnarray}
Now, if $A\ot _RB$ is a twisted tensor product of algebras, then one can
check that the map $T$ given by (\ref{t}) satisfies the axioms in
Theorem \ref{scheme} for $D=A\ot B$, and the deformed multiplication is
the one of $A\ot _RB$, that is $A\ot _RB=(A\ot B)^T$,
so we obtained the associativity of
$A\ot _RB$ as a consequence of Theorem \ref{scheme}.\\
Conversely, if $R:B\otimes A\rightarrow A\otimes B$ is a linear map
such that the map $T$ given by (\ref{t}) is a twistor for $A\otimes
B$, then $R$ is a twisting map and $(A\otimes B)^T=A\otimes _RB$. If
this is the case, we will say that the twistor $T$ is
\dtext{afforded by the twisting map $R$}.

\begin{remark}{\em
If $T$ is a twistor for an algebra $D$, a consequence of (\ref{dec4}) and
(\ref{dec5}) is:
\begin{eqnarray}
&&T(ab\otimes cd)=(a^T)^t(b^{\mathcal{T}})^{\overline{T}}\otimes
(c_\mathcal{T})_T(d_{\overline{T}})_t, \label{bi}
\end{eqnarray}
for all $a, b, c, d\in D$, where $T=t=\mathcal{T}=\overline{T}$.}
\end{remark}

\begin{remark}{\em
Let $T$ be a twistor satisfying the extra conditions
\begin{eqnarray}
&&T_{12}\circ T_{13}=T_{13}\circ T_{12}, \label{cucu1}\\
&&T_{13}\circ T_{23}=T_{23}\circ T_{13}. \label{cucu2}
\end{eqnarray}
Then it is easy to see that $T$ is also an $R$-matrix. Conversely,
a {\it bijective} $R$-matrix satisfying (\ref{cucu1}) and (\ref{cucu2}) is a
twistor. An example of a twistor $T$ satisfying (\ref{cucu1}) and
(\ref{cucu2}) can easily be obtained as follows: take $H$ a cocommutative
bialgebra, $\sigma :H\ot H\rightarrow k$ a bicharacter (i.e. $\sigma $
satisfies $\sigma (1, h)=\sigma (h, 1)=\varepsilon (h)$,
$\sigma (h, h'h'')=\sigma (h_1, h')\sigma (h_2, h'')$ and
$\sigma (hh', h'')=\sigma (h, h''_1)\sigma (h', h''_2)$ for all
$h, h', h''\in H$) and $T:H\ot H\rightarrow H\ot H$,
$T(h\ot h')=\sigma (h_1, h'_1)h_2\ot h'_2$.}
\end{remark}

\begin{remark}{\em
We have seen before (formula (\ref{t})) a basic example of a twistor which
in general is not an $R$-matrix. We present now a basic example of an
$R$-matrix which is not a twistor. Namely, for any algebra $D$, define the
map $T:D\ot D\rightarrow D\ot D$, $T(d\ot d')=d'd\ot 1+1\ot d'd-d'\ot d$.
Then one can check that $T$ is an $R$-matrix (the fact that it satisfies
(\ref{rmat3}) follows from \cite{nuss} or \cite{nichita}) and is 
not a twistor. Note that the multiplication 
$\mu \circ T$ afforded by $T$ is just the multiplication of the opposite
algebra $D^{op}$.}
\end{remark}

\section{More examples of twistors}
\setcounter{equation}{0}
In this section we present more situations where Theorem
\ref{scheme} may be applied.

\noindent\textbf{(i)} Let $A$, $B$, $C$ be three algebras and
$R_1:B\ot A\rightarrow A\ot B$, $R_2:C\ot B\rightarrow B\ot C$,
$R_3:C\ot A\rightarrow A\ot C$ twisting maps. Consider the algebra
$D=A\ot B\ot C$ and the map $T:D\ot D\rightarrow D\ot D$,
\begin{eqnarray}
&&T((a\ot b\ot c)\ot (a'\ot b'\ot c'))=(a\ot b_{R_1}\ot (c_{R_3})_{R_2})\ot
((a'_{R_3})_{R_1}\ot b'_{R_2}\ot c'). \label{itertwistor}
\end{eqnarray}
In general $T$ is {\it not} a twistor for $D$, even if the maps
$R_1$, $R_2$, $R_3$ are compatible. But we have the following result:

\begin{proposition} \label{twisttrei}
With notation as above, $T$ is a twistor for $D$ if and only if the
following conditions hold:
\begin{eqnarray}
&&a_{R_1}\ot (b_{R_1})_{R_2}\ot c_{R_2}=a_{R_1}\ot (b_{R_2})_{R_1}
\ot c_{R_2}, \label{part1} \\
&&(a_{R_1})_{R_3}\ot b_{R_1}\ot c_{R_3}=(a_{R_3})_{R_1}\ot b_{R_1}
\ot c_{R_3}, \label{part2} \\
&&a_{R_3}\ot b_{R_2}\ot (c_{R_3})_{R_2}=a_{R_3}\ot b_{R_2}
\ot (c_{R_2})_{R_3}, \label{part3}
\end{eqnarray}
for all $a\in A$, $b\in B$, $c\in C$. Moreover, in this case it follows
that $R_1$, $R_2$, $R_3$ are compatible twisting maps and
$D^T=A\ot _{R_1}B\ot _{R_2}C$.
\end{proposition}

\begin{proof}
The fact that $T$ is a twistor if and only if (\ref{part1})--(\ref{part3})
hold follows by a direct computation, we leave the details to the reader.
We only prove that $R_1$, $R_2$, $R_3$ are compatible. We compute:
\begin{eqnarray*}
(A\ot R_2)(R_3\ot B)(C\ot R_1)(a\ot b\ot c)&=&
(a_{R_1})_{R_3}\ot (b_{R_1})_{R_2}\ot (c_{R_3})_{R_2}\\
&{{\rm (\ref{part1})}\atop =}&(a_{R_1})_{R_3}\ot (b_{R_2})_{R_1}
\ot (c_{R_3})_{R_2}\\
&{{\rm (\ref{part2})}\atop =}&(a_{R_3})_{R_1}\ot (b_{R_2})_{R_1}
\ot (c_{R_3})_{R_2}\\
&{{\rm (\ref{part3})}\atop =}&(a_{R_3})_{R_1}\ot (b_{R_2})_{R_1}
\ot (c_{R_2})_{R_3}\\
&=&(R_1\ot C)(B\ot R_3)(R_2\ot A)(a\ot b\ot c).
\end{eqnarray*}
The fact that $D^T=A\ot _{R_1}B\ot _{R_2}C$ is obvious.
\end{proof}

\begin{remark} The conditions in Proposition \ref{twisttrei} are
satisfied whenever we start with compatible twisting maps
$R_1$, $R_2$, $R_3$
such that one of them is a usual flip; a concrete example where this
happens is for the so-called two-sided smash product, see \cite{jlpvo}
for details.
\end{remark}

Proposition \ref{twisttrei} may be extended to an iterated twisted
tensor product of any number of factors by means of the Coherence Theorem
stated in \cite{jlpvo}. In order to do this, just realize that conditions
\eqref{part1}, \eqref{part2}, and \eqref{part3} mean simply requiring that
$\{R_1,R_2,\tau _{AC}\}$, $\{R_1,\tau _{BC}, R_3\}$ and
$\{\tau _{AB},R_2,R_3\}$ are
sets of compatible twisting maps, where the $\tau $'s are classical flips.

\begin{proposition}
Let $A_1,\dotsc,A_n$ be some algebras, $\{R_{ij}\}_{i<j}$ a set of twisting
maps,
with $R_{ij}:A_j\ot A_i\to A_i\ot A_j$, and let $D=A_1\ot\dotsb\ot A_n$.
Then the following two conditions are equivalent:
    \begin{enumerate}
        \item The map $T:D\ot D \to D\ot D$ defined by
            \begin{eqnarray*}
                T & := & (Id_{A_1\ot\dotsb\ot A_{n-1}}\ot \tau_{n\,1} \ot
Id_{A_2\ot\dotsb\ot A_{n}}) \circ \dotsb \circ \\
                && \circ (Id_{A_1\ot\dotsb\ot A_{n-k-1}}\ot
\tau_{n-k\,1}\ot\dotsb\ot \tau_{n\, k+1} \ot
Id_{A_{k+2}\ot\dotsb\ot A_{n}})\circ\\
                && \circ\dotsb \circ (Id_{A_1}\ot
\tau_{21}\ot \dotsb\ot \tau _{n\, n-1}\ot Id_{A_n})\circ  (Id_{A_1}\ot
R_{12}\ot \dotsb\ot R _{n-1\, n}\ot Id_{A_n})\circ\\
                && \circ\dotsb\circ (Id_{A_1\ot\dotsb\ot
A_{n-k-1}}\ot R_{1\, n-k}\ot\dotsb\ot R_{k+1\, } \ot
Id_{A_{k+2}\ot\dotsb\ot A_{n}})\circ\dotsb\circ\\
                && \circ (Id_{A_1\ot\dotsb\ot A_{n-1}}\ot
R_{n\,1} \ot Id_{A_2\ot\dotsb\ot A_{n}})
            \end{eqnarray*}
            is a twistor.

        \item For any triple $i<j<k\in\{1,\dotsc,n\}$, we have that
$\{R_{ij},R_{jk},\tau_{ik}\}$, $\{R_{ij},\tau_{jk},R_{ik}\}$ and
$\{\tau_{ij}, R_{jk},R_{ik}\}$ are sets of compatible twisting maps.
    \end{enumerate}
    Moreover, if the conditions are satisfied, then the twisting maps
$\{R_{ij}\}_{i<j}$ are compatible, and we have
$D^T = A_1\ot_{R_{12}}\dotsb\ot _{R_{n-1\, n}}A_n$,
that is, the twisting induced by the twistor $T$ gives the
iterated twisted tensor product associated to the maps.
\end{proposition}

\begin{proof}
        We just outline the main ideas of the proof, leaving
details to the reader. The proof is by induction on the
number of terms $n\geq 3$; for $n=3$, the result is just Proposition
\ref{twisttrei}. Now, assuming the result is true for $n-1$ algebras
with their corresponding twisting maps, and given $A_1,\dotsc,A_n$ algebras,
satisfying the hypothesis of the proposition, we consider the algebras
$B_1:= A_1,\dotsc, B_{n-2}:= A_{n-2},\ B_{n-1}:=
A_{n-1}\ot_{R_{n-1\, n}} A_n$,
with the twisting maps defined as in the Coherence Theorem.
Directly from the hypothesis of the proposition, it follows from the
Coherence Theorem that the newly defined twisting maps also satisfy the
conditions in the proposition, so we may apply our induction
hypothesis to the algebras $B_1,\dotsc, B_{n-1}$.
\end{proof}

A particular case of the former proposition is found in the realization of
the noncommutative planes of Connes and Dubois--Violette as iterated twisted
tensor products (\cite{jlpvo}).
As the twisting maps involved in this process are just
multiples of the classical flips, the compatibility conditions are
trivially satisfied, and the proposition tells us that any noncommutative
$2n$--plane $C_{alg}(\mathbb{R}^{2n}_\theta)$ may also be realized as a
deformation through a twistor of the commutative algebra
$\mathbb{C}[z^1,\bar{z}^1,\dotsc,z^n,\bar{z}^n]$. Moreover, the former
proposition provides an explicit formula for the twistor $T$ that recovers
the iterated twisted tensor product. Taking into account the identification
    \begin{eqnarray*}
        \mathbb{C}[z^1,\bar{z}^1,\dotsc, z^n,\bar{z}^n] &
\longrightarrow & \mathbb{C}[z^1,\bar{z}^1]\ot \dotsb \ot
\mathbb{C}[z^n,\bar{z}^n], \\
        z^i & \longmapsto & 1\ot \dotsb\ot z^i\ot\dotsb\ot 1, \\
        \bar{z}^i & \longmapsto & 1\ot \dotsb\ot
\bar{z}^i\ot\dotsb\ot 1,
    \end{eqnarray*}
    where $z^i$ maps to the position $2i-1$ and $\bar{z}^i$ maps to
the position $2i$, it is easy to realize that the twistor given by the
proposition is defined on generators as:
    \begin{align*}
        T(z^i\ot z^j) & =
            \begin{cases}
                z^i\ot z^j & \text{if $i\leq j$,}\\
                \lambda^ {ij} z^i\ot z^j & \text{otherwise},
            \end{cases} &
        T(\bar{z}^i\ot \bar{z}^j) & =
            \begin{cases}
                \bar{z}^i\ot \bar{z}^j &
\text{if $i\leq j$,}\\
                \lambda^ {ij} \bar{z}^i\ot \bar{z}^j &
\text{otherwise},
            \end{cases} \\
        T(\bar{z}^i\ot z^j) & =
            \begin{cases}
                \bar{z}^i\ot z^j & \text{if $i\leq j$,}\\
                \lambda^{ji} \bar{z}^i\ot z^j &
\text{otherwise},
            \end{cases} &
        T(z^i\ot \bar{z}^j) & =
            \begin{cases}
                z^i\ot \bar{z}^j & \text{if $i\leq j$,}\\
                \lambda^ {ji} z^i\ot \bar{z}^j &
\text{otherwise}.
            \end{cases}
    \end{align*}

\noindent\textbf{(ii)} Let $A$ be an algebra with multiplication
$\mu _A=\mu $ and $H$ a bialgebra such that $A$ is an $H$-bimodule
algebra with actions denoted by $\pi _l:H\ot A\rightarrow A$, $\pi
_l(h\ot a)=h\cdot a$ and $\pi _r:A\ot H\rightarrow A$, $\pi _r(a\ot
h)=a\cdot h$, also $A$ is an $H$-bicomodule algebra, with coactions
denoted by $\psi _l:A\rightarrow H\ot A$, $a\mapsto a_{[-1]}\ot
a_{[0]}$ and $\psi _r:A\rightarrow A\ot H$, $a\mapsto a_{<0>}\ot
a_{<1>}$, and moreover the following compatibility conditions hold,
for all $h\in H$ and $a\in A$:
\begin{eqnarray*}
&&(h\cdot a)_{[-1]}\ot (h\cdot a)_{[0]}=a_{[-1]}\ot h\cdot a_{[0]}, \;\;\;
(h\cdot a)_{<0>}\ot (h\cdot a)_{<1>}=h\cdot a_{<0>}\ot a_{<1>}, \\
&&(a\cdot h)_{[-1]}\ot (a\cdot h)_{[0]}=a_{[-1]}\ot a_{[0]}\cdot h, \;\;\;
(a\cdot h)_{<0>}\ot (a\cdot h)_{<1>}=a_{<0>}\cdot h\ot a_{<1>}.
\end{eqnarray*}
Such a datum was considered in \cite{pvo}, where it is called
an L-R-twisting datum for $A$ (and contains as particular case the concept
of very strong left twisting datum from \cite{santos}, which is obtained
if the right action and coaction are trivial).
\begin{proposition} (\cite{pvo})  \label{lrtwist}
Given an L-R-twisting datum, define a new multiplication
on $A$ by
\begin{eqnarray}
&&a\bullet a'=(a_{[0]}\cdot a'_{<1>})(a_{[-1]}\cdot a'_{<0>}), \;\;\;
\forall \;\;a, a'\in A. \label{bullet}
\end{eqnarray}
Then $(A, \bullet , 1)$ is an associative unital algebra.
\end{proposition}
This result may be obtained as a consequence of Theorem \ref{scheme}.
Namely, define
\begin{eqnarray}
&&T:A\ot A\rightarrow A\ot A, \;\;\;T(a\ot a')=a_{[0]}\cdot a'_{<1>}\ot
a_{[-1]}\cdot a'_{<0>}. \label{twlr}
\end{eqnarray}
Then one can check that $T$ is a twistor for $A$, and obviously the
new multiplication $\bullet $ defined above coincides with $\mu
\circ T$.

\noindent\textbf{(iii)} Let $H$, $K$ be two bialgebras, $A$ an
algebra which is a left $H$-comodule algebra with coaction $a\mapsto
a_{[-1]}\ot a_{[0]}\in H\ot A$ and a left $K$-module algebra with
action $k\ot a\mapsto k\cdot a$, for all $a\in A$, $k\in K$, such
that $(k\cdot a)_{[-1]}\ot (k\cdot a)_{[0]}=a_{[-1]}\ot k\cdot
a_{[0]}$, for all $a\in A$, $k\in K$. Let $f:H\rightarrow K$ be a
bialgebra map. Then, by \cite{cj}, the new multiplication defined on
$A$ by $a\cdot _fa'=a_{[0]}(f(a_{[-1]}) \cdot a')$ is associative
with unit $1$. This multiplication is afforded by the map $T:A\ot
A\rightarrow A\ot A$, $T(a\ot a')=a_{[0]}\ot f(a_{[-1]})\cdot a'$,
which is easily seen to be a twistor.

\noindent\textbf{(iv)} Let $H$ be a bialgebra and $F=F^1\ot F^2\in
H\ot H$ an element with $(\varepsilon \ot H)(F)=(H\ot \varepsilon
)(F)=1$. Assume that $F$ satisfies the following list of axioms,
considered in \cite{jc}, \cite{km}: $(H\ot \Delta
)(F)=F_{13}F_{12}$, $(\Delta \ot H)(F)=F_{13}F_{23}$ and
$F_{12}F_{23}=F_{23}F_{12}$. Let $D$ be a left $H$-module algebra
and define $T:D\ot D\rightarrow D\ot D$ by $T(d\ot d')=F^1\cdot d\ot
F^2\cdot d'$. Then it is easy to see that $T$ is a twistor for $D$.
In case $F$ is invertible, the multiplication of $D^T$ fits into the
well-known procedure of twisting a module algebra by a Drinfeld
twist.

\noindent\textbf{(v)} Let $H$ be a bialgebra and $\sigma :H\otimes
H\rightarrow k$ a linear map. Define $T:H\otimes H\rightarrow
H\otimes H$ by $T(a\otimes b)=\sigma (a_1, b_1)a_2\otimes b_2$, for
all $a, b\in H$. Then $T$ is a twistor for $H$ if and only if
$\sigma $ satisfies the following conditions: $\sigma (a,
1)=\varepsilon (a)=\sigma (1, a)$, $\sigma (a, bc)=\sigma (a_1,
b)\sigma (a_2, c)$, $\sigma (ab, c)=\sigma (a, c_2)\sigma (b, c_1)$
and $\sigma (a, b_1)\sigma (b_2, c)=\sigma (b_1, c)\sigma (a, b_2)$,
for all $a, b, c\in H$. Note that elements satisfying the last
condition have been considered in \cite{psvo}, under the name {\it
neat} elements.

\noindent\textbf{(vi)} Let $(D, \delta )$ be a differential
associative algebra, that is $D$ is an associative algebra and
$\delta :D\rightarrow D$ is a derivation (i.e. $\delta (dd')=\delta
(d)d'+d\delta (d')$) with $\delta ^2=0$. Then one can see that the
map $T:D\otimes D\rightarrow D\otimes D$, $T(d\otimes d')=d\otimes
d'+\delta (d)\otimes \delta (d')$ is a twistor for $D$.

\section{Some properties of twistors}
\setcounter{equation}{0}
\begin{proposition}\label{morfism}
Let $T$ be a twistor for an algebra $D$ and $U$ a twistor for an algebra $F$.
If $\nu :D\rightarrow F$ is an algebra map such that $(\nu \otimes \nu )
\circ T=U\circ (\nu \otimes \nu )$, then $\nu $ is also an algebra map
from $D^T$ to $F^U$.
\end{proposition}

It was proved in \cite{bor} that, if $A\otimes _R B$ and $A'\otimes
_{R'} B'$ are twisted tensor products of algebras and
$f:A\rightarrow A'$ and $g:B\rightarrow B'$ are algebra maps
satisfying the condition $(f\otimes g)\circ R=R'\circ (g\otimes f)$,
then $f\otimes g:A\otimes _R B\rightarrow A'\otimes _{R'} B'$ is an
algebra map. One can easily see that this result is a particular
case of Proposition \ref{morfism}, with $D=A\otimes B$, $F=A'\otimes
B'$, $\nu =f\otimes g$ and $T$ (respectively $U$) the twistor
afforded by $R$ (respectively $R'$).

 We present one more situation where Proposition \ref{morfism} may 
be applied. We recall that the 
L-R-smash product over a cocommutative Hopf algebra was introduced
in \cite{b2}, \cite{b3}, and generalized to an arbitrary Hopf
algebra in \cite{pvo} as follows: if $A$ is an $H$-bimodule algebra,
the L-R-smash product $A\nat H$ is the following algebra structure
on $A\ot H$:
\begin{eqnarray*}
&&(a\nat h)(a'\nat h')=(a\cdot h'_2)(h_1\cdot a')
\nat h_2h'_1, \;\;\;\forall \;a, a'\in A, \;
h, h'\in H.
\end{eqnarray*}
The diagonal crossed product $A\bowtie H$ is the following algebra
structure on $A\ot H$, see \cite{hn}, \cite{bpvo}:
\begin{eqnarray*}
&&(a\bowtie h)(a'\bowtie h')=
a(h_1\cdot a'\cdot S^{-1}(h_3))
\bowtie h_2h', \;\;\;\forall \;a, a'\in A, \;
h, h'\in H.
\end{eqnarray*}

It was proved in \cite{pvo} that actually $A\bowtie H$ and $A\nat H$
are isomorphic as algebras. This result may be reobtained using
Proposition \ref{morfism} as follows. Denote by $A\# _rH$ the
algebra structure on $A\otimes H$ with multiplication $(a\otimes
h)(a'\otimes h')=(a\cdot h'_2)a'\otimes hh'_1$, and by $A\bowtie
_rH$ the algebra structure on $A\otimes H$ with multiplication
$(a\otimes h)(a'\otimes h')=a(a'\cdot S^{-1}(h_2))\otimes h_1h'$.
One may check that the map $\nu :A\bowtie _rH\rightarrow A\# _rH$
given by $\nu (a\otimes h)=a\cdot h_2\otimes h_1$ is an algebra map
(actually, an isomorphism, with inverse $\nu ^{-1}(a\otimes
h)=a\cdot S^{-1}(h_2) \otimes h_1$). Define now the map $T:(A\otimes
H)\otimes (A\otimes H) \rightarrow (A\otimes H)\otimes (A\otimes H)$
by $T((a\otimes h)\otimes (a'\otimes h'))=(a\otimes h_2)\otimes
(h_1\cdot a'\otimes h')$. Then one may check, by direct computation,
that $T$ is a twistor for both $A\# _rH$ and $A\bowtie _rH$, and
moreover $(A\# _rH)^T=A\nat H$, $(A\bowtie _rH)^T=A\bowtie H$ and
$(\nu \otimes \nu )\circ T= T\circ (\nu \otimes \nu )$. Hence,
Proposition \ref{morfism} may be applied and we obtain as a
consequence that $\nu $ is an algebra map from $A\bowtie H$ to
$A\nat H$.

By \cite{pvo}, the L-R-twisted product (\ref{bullet})
may be obtained as a left twisting followed by a right twisting and
viceversa. This fact admits an interpretation in terms
of twistors.

\begin{proposition} \label{stdr}
Let $D$ be an algebra and $X, Y:D\ot D\rightarrow D\ot D$ two twistors for
$D$, satisfying the following conditions:
\begin{eqnarray}
&&X_{23}\circ Y_{12}=Y_{12}\circ X_{23}, \label{l1} \\
&&X_{23}\circ Y_{13}=Y_{13}\circ X_{23}, \label{l2} \\
&&X_{12}\circ Y_{23}=Y_{23}\circ X_{12}, \label{l3} \\
&&X_{12}\circ Y_{13}=Y_{13}\circ X_{12}. \label{l4}
\end{eqnarray}
Then $Y$ is a twistor for $D^X$, $X$ is a twistor for $D^Y$, $X\circ Y$ and
$Y\circ X$ are twistors for $D$ and of course $(D^X)^Y=D^{X\circ Y}$ and
$(D^Y)^X=D^{Y\circ X}$.
\end{proposition}

\begin{proof}
Note first that (\ref{l2}) and (\ref{l4}) are respectively equivalent to
$X_{13}\circ Y_{23}=Y_{23}\circ X_{13}$ and $Y_{12}\circ X_{13}=
X_{13}\circ Y_{12}$, hence the above conditions are actually symmetric in
$X$ and $Y$, so we only have to prove that $Y$ is a twistor for $D^X$ and
$X\circ Y$ is a twistor for $D$. \\
To prove that $Y$ is a twistor for $D^X$ we only have to check
(\ref{dec4}) and (\ref{dec5}) for $Y$ with respect to the multiplication
$*$ of $D^X$; we compute:
\begin{eqnarray*}
d^Y\ot (d'*d'')_Y&=&d^Y\ot (d'^Xd''_X)_Y\\
&{{\rm (\ref{dec4})}\atop =}&(d^Y)^y\ot (d'^X)_Y(d''_X)_y\\
&{{\rm (\ref{l2})}\atop =}&(d^Y)^y\ot (d'^X)_Y(d''_y)_X\\
&{{\rm (\ref{l1})}\atop =}&(d^Y)^y\ot (d'_Y)^X(d''_y)_X\\
&=&(d^Y)^y\ot d'_Y*d''_y,
\end{eqnarray*}
\begin{eqnarray*}
(d*d')^Y\ot d''_Y&=&(d^Xd'_X)^Y\ot d''_Y\\
&{{\rm (\ref{dec5})}\atop =}&(d^X)^Y(d'_X)^y\ot (d''_y)_Y\\
&{{\rm (\ref{l3})}\atop =}&(d^X)^Y(d'^y)_X\ot (d''_y)_Y\\
&{{\rm (\ref{l4})}\atop =}&(d^Y)^X(d'^y)_X\ot (d''_y)_Y\\
&=&d^Y*d'^y\ot (d''_y)_Y.
\end{eqnarray*}
Now we check (\ref{dec4}) and (\ref{dec5}) for $T:=X\circ Y$; we compute:
\begin{eqnarray*}
d^T\ot (d'd'')_T&=&(d^Y)^X\ot ((d'd'')_Y)_X\\
&{{\rm (\ref{dec4})}\atop =}&((d^Y)^y)^X\ot (d'_Yd''_y)_X\\
&{{\rm (\ref{dec4})}\atop =}&(((d^Y)^y)^X)^x\ot (d'_Y)_X(d''_y)_x\\
&{{\rm (\ref{l4})}\atop =}&(((d^Y)^X)^y)^x\ot (d'_Y)_X(d''_y)_x\\
&=&(d^T)^t\ot d'_Td''_t,
\end{eqnarray*}
\begin{eqnarray*}
(dd')^T\ot d''_T&=&((dd')^Y)^X\ot (d''_Y)_X\\
&{{\rm (\ref{dec5})}\atop =}&(d^Yd'^y)^X\ot ((d''_y)_Y)_X\\
&{{\rm (\ref{dec5})}\atop =}&(d^Y)^X(d'^y)^x\ot (((d''_y)_Y)_x)_X\\
&{{\rm (\ref{l2})}\atop =}&(d^Y)^X(d'^y)^x\ot (((d''_y)_x)_Y)_X\\
&=&d^Td'^t\ot (d''_t)_T.
\end{eqnarray*}
It remains to prove (\ref{dec3}) for $T$; we compute:
\begin{eqnarray*}
T_{12}\circ T_{23}&=&X_{12}\circ Y_{12}\circ X_{23}\circ Y_{23}\\
&{{\rm (\ref{l1})}\atop =}&X_{12}\circ X_{23}\circ Y_{12}\circ Y_{23}\\
&{{\rm (\ref{dec3})}\atop =}&X_{23}\circ X_{12}\circ Y_{23}\circ Y_{12}\\
&{{\rm (\ref{l3})}\atop =}&X_{23}\circ Y_{23}\circ X_{12}\circ Y_{12}\\
&=&T_{23}\circ T_{12},
\end{eqnarray*}
and the proof is finished.
\end{proof}

Let now $A$ be as in Proposition \ref{lrtwist} and define
$X, Y:A\ot A\rightarrow A\ot A$ by
\begin{eqnarray*}
&&X(a\ot a')=a\cdot a'_{<1>}\ot a'_{<0>}, \;\;\;
Y(a\ot a')=a_{[0]}\ot a_{[-1]}\cdot a'.
\end{eqnarray*}
Then one can check that $X$ and $Y$ satisfy the hypotheses of
Proposition \ref{stdr}, and moreover we have $X\circ Y=Y\circ X=T$,
where $T$ is given by (\ref{twlr}). Hence, we obtain $(A, \bullet ,
1)=(A^X)^Y=(A^Y)^X$.

 Also as a consequence of Proposition
\ref{stdr}, we obtain that if $T$ is a twistor for an algebra $D$,
satisfying (\ref{cucu1}) and (\ref{cucu2}), then $T$ is a twistor
also for $D^T$, hence we obtain a sequence of associative algebras
$D$, $D^T$, $D^{T^2}$, $D^{T^3}$, etc.

 A particular case of
Proposition \ref{stdr} is the following:
\begin{corollary}
 Let $A$, $B$ be two algebras and $R, S:B\otimes A\rightarrow A\otimes B$
two twisting maps. Denote by $X$ (respectively $Y$) the twistor for
$A\otimes B$ afforded by $R$ (respectively $S$) and assume that the
following conditions are satisfied:
\begin{eqnarray*}
&&(a_R)_S\otimes b_R\otimes b'_S=(a_S)_R\otimes b_R\otimes b'_S, \;\;\;
a_R\otimes a'_S\otimes (b_R)_S=a_R\otimes a'_S\otimes (b_S)_R,
\end{eqnarray*}
for all $a, a'\in A$ and $b, b'\in B$. Define $R*S, S*R:B\otimes A
\rightarrow A\otimes B$ by
\begin{eqnarray*}
(R*S)(b\otimes a)=(a_S)_R\otimes (b_S)_R, \;\;\;
(S*R)(b\otimes a)=(a_R)_S\otimes (b_R)_S.
\end{eqnarray*}
Then $Y$ is a twistor for $A\otimes _R B$, $X$ is a twistor for
$A\otimes _S B$, $X\circ Y$ (respectively $Y\circ X$) is a twistor for
$A\otimes B$ afforded by the twisting map $R*S$ (respectively $S*R$) and
we have $(A\otimes _R B)^Y=A\otimes _{R*S} B$,
$(A\otimes _S B)^X=A\otimes _{S*R} B$.
\end{corollary}

We are now interested in lifting
(bi) module structures from an algebra $D$ to $D^T$.
This is achieved in the next
result, the proof follows from a direct computation and will be omitted.

\begin{proposition} \label{modtw}
Let $D$ be an algebra and $T$ a twistor for $D$.

\noindent\textbf{(i)} Let $V$ be a left $D$-module, with action
$\lambda :D\ot V\rightarrow V$, $\lambda (d\ot v)=d\cdot v$. Assume
that we are given a linear map $\Gamma :D\ot V\rightarrow D\ot V$,
with notation $\Gamma (d\ot v)=d_{\Gamma }\ot v_{\Gamma }$, for all
$d\in D$, $v\in V$, such that $\Gamma (1\ot v)=1\ot v$, for all
$v\in V$, and
\begin{eqnarray}
&&\lambda _{23}\circ \Gamma _{13}\circ T_{12}=\Gamma \circ \lambda _{23}
:D\ot D\ot V\rightarrow D\ot V, \label{m1} \\
&&\mu _{12}\circ \Gamma _{13}\circ \Gamma _{23}=\Gamma \circ \mu _{12}
:D\ot D\ot V\rightarrow D\ot V, \label{m2} \\
&&T_{12}\circ \Gamma _{23}=\Gamma _{23}\circ T_{12}:D\ot D\ot V
\rightarrow D\ot D\ot V. \label{m3}
\end{eqnarray}
Then $V$ becomes a left $D^T$-module, with action $\lambda \circ
\Gamma : D\ot V\rightarrow V$. We denote by $V^{\Gamma }$ this
$D^T$-module structure on $V$ and by $d\rightarrow v=d_{\Gamma
}\cdot v_{\Gamma }$ the action of $D^T$ on $V$. We call the map
$\Gamma $ a \dtext{left module twistor} for $V$ relative to $T$.

\noindent\textbf{(ii)} Let $V$ be a right $D$-module, with action
$\rho :V\otimes D\rightarrow V$, $\rho (v\otimes d)=v\cdot d$, and
assume that we are given a linear map $\Pi :V\otimes D\rightarrow
V\otimes D$, with notation $\Pi (v\otimes d)=v_{\Pi }\otimes d_{\Pi
}$, for all $d\in D$, $v\in V$, such that $\Pi (v\otimes 1)=v\otimes
1$, for all $v\in V$, and
\begin{eqnarray}
&&\mu _{23}\circ \Pi _{13}\circ \Pi _{12}=\Pi \circ \mu _{23}:V\otimes D
\otimes D\rightarrow V\otimes D, \label{n1} \\
&&\rho _{12}\circ \Pi _{13}\circ T_{23}=\Pi \circ \rho _{12}:V\otimes D
\otimes D\rightarrow V\otimes D, \label{n2}  \\
&&\Pi _{12}\circ T_{23}=T_{23}\circ \Pi _{12}:V\otimes D\otimes D
\rightarrow V\otimes D\otimes D. \label{n3}
\end{eqnarray}
Then $V$ becomes a right $D^T$-module, with action $\rho \circ \Pi :
V\otimes D\rightarrow V$. We denote by $^{\Pi }V$ this $D^T$-module
structure on $V$ and by $v\leftarrow d=v_{\Pi }\cdot d_{\Pi }$ the
action of $D^T$ on $V$. We call the map $\Pi $ a \dtext{right module
twistor} for $V$ relative to $T$.

\noindent\textbf{(iii)} Let $V$ be a $D$-bimodule, and let $\Gamma $
and $\Pi $ be a left respectively a right module twistor for $V$
relative to $T$. Assume that the following conditions hold:
\begin{eqnarray}
&&\rho _{23}\circ T_{13}\circ \Gamma _{12}=\Gamma \circ \rho _{23}
:D\otimes V\otimes D\rightarrow D\otimes V, \label{bim1} \\
&&\lambda _{12}\circ T_{13}\circ \Pi _{23}=\Pi \circ \lambda _{12}
:D\otimes V\otimes D\rightarrow V\otimes D, \label{bim2} \\
&&\Gamma _{12}\circ \Pi _{23}=\Pi _{23}\circ \Gamma _{12}
:D\otimes V\otimes D\rightarrow D\otimes V\otimes D. \label{bim3}
\end{eqnarray}
Let $^{\Pi }V^{\Gamma }$ be $V^{\Gamma }$ as a left $D^T$-module and
$^{\Pi }V$ as a right $D^T$-module. Then $^{\Pi }V^{\Gamma }$ is a
$D^T$-bimodule.
\end{proposition}

We recall from \cite{cap} the following result. Let $A\ot _RB$ be a twisted
tensor product of algebras, $M$ a left $A$-module, $N$ a left $B$-module
(we denote by $\lambda _M$ and respectively $\lambda _N$ the actions) and
$\tau _{M, B}:B\ot M\rightarrow M\ot B$ a linear map, with notation
$\tau _{M, B}(b\ot m)=m_{\tau }\ot b_{\tau }$, such that
$\tau _{M, B}(1\ot m)=m\ot 1$, for all $m\in M$, and the following
conditions hold:
\begin{eqnarray*}
&&\tau _{M, B}\circ (\mu _B\ot M)=(M\ot \mu _B)\circ (\tau _{M, B}\ot B)
\circ (B\ot \tau _{M, B}), \\
&&\tau _{M, B}\circ (B\ot \lambda _M)=(\lambda _M\ot B)\circ (A\ot
\tau _{M, B})\circ (R\ot M)
\end{eqnarray*}
(such a map $\tau _{M, B}$ is called a left module twisting map).
Then $M\ot N$ becomes a left $A\ot _RB$-module, with action $(a\ot
b)\rightarrow (m\ot n)=a\cdot m_{\tau }\ot b_{\tau }\cdot n$. This
result is a particular case of Proposition \ref{modtw} (i). Indeed,
we consider the algebra $D=A\ot B$ (the ordinary tensor product),
the twistor $T$ for $D$ given by (\ref{t}), the left $D$-module
$V=M\ot N$ with action $(a\ot b)\cdot (m\ot n)=a\cdot m\ot b\cdot
n$, and the map $\Gamma :(A\ot B)\ot (M\ot N)\rightarrow (A\ot B)\ot
(M\ot N)$ given by $\Gamma ((a\ot b)\ot (m\ot n))=(a\ot b_{\tau
})\ot (m_{\tau }\ot n)$. Then one can check that $\Gamma $ satifies
the axioms of a left module twistor, and the left $D^T=A\ot
_RB$-module $V^{\Gamma }$ is obviously the $A\otimes _R B$-module
structure on $M\ot N$ presented above. Similarly, one can see that
Proposition \ref{modtw} (ii) contains as particular case the lifting
of right module structures to a twisted tensor product from
\cite{cap}.

Another example may be obtained as follows. Let $A$ be as in
Proposition \ref{lrtwist}, and $V$ a vector space which is a left
$A$-module (with action $a\ot v\mapsto a\cdot v$), a left $H$-module
(with action $h\ot v\mapsto h\cdot v$) and a right $H$-comodule
(with coaction $v\mapsto v_{<0>}\ot v_{<1>}\in V\ot H$) such that
the following conditions are satisfied, for all $h\in H$, $a\in A$,
$v\in V$:
\begin{eqnarray*}
&&(h\cdot v)_{<0>}\ot (h\cdot v)_{<1>}=h\cdot v_{<0>}\ot v_{<1>}, \\
&& h\cdot (a\cdot v)=(h_1\cdot a)\cdot (h_2\cdot v), \\
&&(a\cdot v)_{<0>}\ot (a\cdot v)_{<1>}=a_{<0>}\cdot v_{<0>}\ot
a_{<1>}v_{<1>}.
\end{eqnarray*}
Define the map $\Gamma :A\ot V\rightarrow A\ot V$ by $\Gamma (a\ot
v)=a_{[0]}\cdot v_{<1>}\ot a_{[-1]}\cdot v_{<0>}$. Then one can
check that $\Gamma $ and the twistor $T$ given by (\ref{twlr})
satisfy the hypotheses of Proposition \ref{modtw} (i), hence $V$
becomes a left module over $(A, \bullet )$, with action
$a\rightarrow v=(a_{[0]}\cdot v_{<1>})\cdot (a_{[-1]}\cdot
v_{<0>})$.

 We present now an application of Proposition \ref{modtw}.
\begin{proposition}
Let $(D, \mu , u)$ be an algebra and consider the universal first order
differential calculus $\Omega _u^1(D)=Ker (\mu )$, with its canonical
$D$-bimodule structure. If $T$ is a twistor for $D$, then $\Omega _u^1(D)$
becomes also a $D^T$-bimodule.
\end{proposition}

\begin{proof}
Define the maps $\Gamma , \Pi :D\otimes D\otimes D\rightarrow
D\otimes D\otimes D$ by $\Gamma =T_{13}\circ T_{12}$ and $\Pi =T_{13}\circ
T_{23}$. We claim that $\Gamma (D\otimes Ker (\mu ))\subseteq
D\otimes Ker (\mu )$ and $\Pi (Ker (\mu )\otimes D)\subseteq Ker (\mu )
\otimes D$. To prove this, we recall the following result from linear
algebra: if $f:V\rightarrow V'$ and $g:W\rightarrow W'$ are linear
maps, then $Ker (f\otimes g)=Ker (f)\otimes W+V\otimes Ker (g)$. We apply
this result for the map $D\otimes \mu :D\otimes D\otimes D\rightarrow
D\otimes D\otimes D$, and we obtain $Ker (D\otimes \mu )=Ker (D)\otimes D
\otimes D+D\otimes Ker (\mu )=D\otimes Ker (\mu )$. Let $x\in D\otimes
Ker (\mu )$; in order to prove that $\Gamma (x)\in D\otimes Ker (\mu )$,
in view of the above it is enough to prove that $((D\otimes \mu )\circ
\Gamma )(x)=0$. But using (\ref{dec1}) and the definition of $\Gamma $,
we see that $(D\otimes \mu )\circ \Gamma =T\circ \mu _{23}$, and
obviously $(T\circ \mu _{23})(x)=0$ because $x\in D\otimes Ker (\mu )$.
Similarly one can prove that $\Pi (Ker (\mu )\otimes D)\subseteq
Ker (\mu )\otimes D$. Now, if we denote by $\lambda :D\otimes Ker (\mu )
\rightarrow Ker (\mu )$ and $\rho :Ker (\mu )\otimes D\rightarrow
Ker (\mu )$ the left and right $D$-module structures of $Ker (\mu )$
(given by $\lambda =\mu _{12}$ and $\rho =\mu _{23}$), then the
maps $\lambda , \rho , \Gamma , \Pi $ satisfy all the hypotheses of
Proposition \ref{modtw} (this proof is a direct computation and is
omitted), hence indeed $Ker (\mu )$ becomes a $D^T$-bimodule.
\end{proof}

Actually, more can be said about this $D^T$-bimodule $Ker (\mu )$.
Denote by $\delta :D\rightarrow Ker (\mu )$, $\delta (d)=d\otimes
1-1\otimes d$ the canonical $D$-derivation.

\begin{proposition}
This map $\delta $ is also a $D^T$-derivation from $D^T$ to $Ker (\mu )$,
where the $D^T$-bimodule structure on $Ker (\mu )$ is the one
presented above.
\end{proposition}

\begin{proof}
Using the formulae for $\Gamma $ and $\Pi $, one can easily see
that $d\rightarrow \delta (d')=d^T\cdot \delta (d'_T)$ and
$\delta (d)\leftarrow d'=\delta (d^T)\cdot d'_T$ for all $d, d'\in D$,
so we immediately obtain:
\begin{eqnarray*}
\delta (d*d')&=&\delta (d^Td'_T)\\
&=&d^T\cdot \delta (d'_T)+\delta (d^T)\cdot d'_T\\
&=&d\rightarrow \delta (d')+\delta (d)\leftarrow d',
\end{eqnarray*}
finishing the proof.
\end{proof}

\begin{proposition}
If the twistor $T$ is bijective, then $(Ker (\mu ), \delta )$ is also a
first order differential calculus over the algebra $D^T$.
\end{proposition}

\begin{proof}
We only have to prove that $Ker (\mu )$ is generated by
$\{\delta (d): d\in D\}$ as a $D^T$-bimodule. If $d, d'\in D$, we denote by
$T^{-1}(d\otimes d')=d^U\otimes d'_U$. If $x=\sum _ia_i\otimes b_i\in
Ker (\mu )$, we can write $x=\sum _i \delta (a_i)\cdot b_i$, which in turn
may be written as $x=\sum _i \delta (a_i^U)\leftarrow (b_i)_U$, q.e.d.
\end{proof}

\section{Pseudotwistors and braided (graded) twistors}
\label{secbraided}
\setcounter{equation}{0}
Let $(\Omega , d)$ be a DG algebra, that is $\Omega =\bigoplus _{n\geq 0}
\Omega ^n$ is a graded algebra and $d:\Omega \rightarrow \Omega $ is a
linear map with $d(\Omega ^n)\subseteq \Omega ^{n+1}$ for all
$n\geq 0$, $d^2=0$ and $d(\omega \zeta )=d(\omega )\zeta +(-1)^{|\omega |}
\omega d(\zeta )$ for all homogeneous $\omega $ and $\zeta $, where
$|\omega |$ is the degree of $\omega $. The Fedosov product
(\cite{fedosov}, \cite{cuqui}), given by
\begin{eqnarray}
&&\omega \circ \zeta =\omega \zeta -(-1)^{|\omega |}d(\omega )d(\zeta ),
\label{fed}
\end{eqnarray}
for homogeneous $\omega $ and $\zeta $, defines a new associative algebra
structure on $\Omega $. If we define the map
\begin{eqnarray}
&&T:\Omega \otimes \Omega \rightarrow \Omega \otimes \Omega , \;\;\;
T(\omega \otimes \zeta )=\omega \otimes \zeta -(-1)^{|\omega |}
d(\omega )\otimes d(\zeta ), \label{fedtw}
\end{eqnarray}
then $T$ satisfies (\ref{dec3}) but fails to satisfy (\ref{dec1}) and
(\ref{dec2}). However, the failure is only caused by some signs, so we
were led to introduce a graded analogue of a twistor, which in turn
leads us to the following much more general concept:
\begin{proposition} \label{tilda}
Let $\mathcal{C}$ be a (strict) monoidal category, $A$ an algebra in
$\mathcal{C}$ with multiplication $\mu $ and unit $u$, $T:A\otimes
A\rightarrow A\otimes A$ a morphism in $\mathcal{C}$ such that
$T\circ (u \otimes A)=u\otimes A$ and $T\circ (A\otimes u)=A\otimes
u$. Assume that there exist two morphisms $\tilde{T}_1,
\tilde{T}_2:A\otimes A\otimes A \rightarrow A\otimes A\otimes A$ in
$\mathcal{C}$ such that
\begin{eqnarray}
&&(A\otimes \mu )\circ \tilde{T}_1\circ (T\otimes A)=T\circ (A\otimes \mu ),
\label{catw1} \\
&&(\mu \otimes A)\circ \tilde{T}_2\circ (A\otimes T)=T\circ (\mu \otimes A),
\label{catw2} \\
&&\tilde{T}_1\circ (T\otimes A)\circ (A\otimes T)=
\tilde{T}_2\circ (A\otimes T)\circ (T\otimes A). \label{catw3}
\end{eqnarray}
Then $(A, \mu \circ T, u)$ is also an algebra in ${\mathcal{C}}$,
denoted by $A^T$. The morphism $T$ is called a
\dtext{pseudotwistor} and the two morphisms $\tilde{T}_1$,
$\tilde{T}_2$ are called the \dtext{companions} of $T$.
\end{proposition}
\begin{proof}
Obviously $u$ is a unit for $(A, \mu \circ T)$, so we
only check the associativity of $\mu \circ T$:
\begin{eqnarray*}
(\mu \circ T)\circ ((\mu \circ T)\otimes A)&=&(\mu \circ T)\circ
(\mu \otimes A)\circ (T\otimes A)\\
&{{\rm (\ref{catw2})}\atop =}&\mu \circ (\mu \otimes A)\circ \tilde{T}_2\circ 
(A\otimes T)\circ (T\otimes A)\\
&{{\rm (\ref{catw3})}\atop =}&\mu \circ (\mu \otimes A)\circ \tilde{T}_1\circ
(T\otimes A)\circ (A\otimes T)\\
&=&\mu \circ (A\otimes \mu )\circ \tilde{T}_1\circ (T\otimes A)\circ
(A\otimes T)\\
&{{\rm (\ref{catw1})}\atop =}&\mu \circ T\circ (A\otimes \mu )
\circ (A\otimes T)\\
&=&(\mu \circ T)\circ (A\otimes (\mu \circ T)),
\end{eqnarray*}
finishing the proof.
\end{proof}

\begin{remark}{\em
Obviously, an ordinary twistor $T$ is a pseudotwistor with companions
$\tilde{T}_1=\tilde{T}_2=T_{13}$. Also, if
$T:A\otimes A\rightarrow A\otimes A$ is a {\it bijective} R-matrix,
one can easily check that $T$ is a pseudotwistor, with companions
$\tilde{T}_1=T_{12}\circ T_{13}\circ T_{12}^{-1}$ and
$\tilde{T}_2=T_{23}\circ T_{13}\circ T_{23}^{-1}$.}
\end{remark}

A pseudotwistor may be thought of as some sort of analogue of a (Hopf)
2-cocycle, as suggested by the following examples (for which $\mathcal{C}$
is the usual category of vector spaces):
\begin{examples} \label{drinfeld} {\em
Let $H$ be a bialgebra and $F=F^1\otimes F^2=f^1\otimes f^2\in H\otimes H$ 
a Drinfeld 
twist, i.e. an invertible element (with inverse denoted by $F^{-1}=
G^1\otimes G^2$) such that $F^1f^1_1\otimes F^2f^1_2\otimes f^2=
f^1\otimes F^1f^2_1\otimes F^2f^2_2$ and 
$(\varepsilon \otimes H)(F)=(H\otimes 
\varepsilon )(F)=1$. If $A$ is a left $H$-module algebra, it is well-known
that the new product on $A$ given by $a*b=(G^1\cdot a)(G^2\cdot b)$ is
associative. This product is afforded by the map $T:A\otimes A
\rightarrow A\otimes A$, $T(a\otimes b)=G^1\cdot a\otimes G^2\cdot b$,
and one may check that $T$ is a pseudotwistor with companions
$\tilde{T}_1$, $\tilde{T}_2$ given by the formulae
\begin{eqnarray*}
&&\tilde{T}_1(a\otimes
b\otimes c)=G^1F^1\cdot a\otimes G^2_1F^2\cdot b\otimes G^2_2\cdot c, \\
&&\tilde{T}_2(a\otimes b\otimes c)=G^1_1\cdot a\otimes G^1_2F^1\cdot b
\otimes G^2F^2\cdot c.
\end{eqnarray*}
Dually, let $H$ be a bialgebra and $\sigma :H\otimes H\rightarrow k$ a
normalized and convolution invertible left 2-cocycle (i.e.
$\sigma $ satisfies $\sigma (h_1, h'_1)\sigma (h_2h'_2, h'')=
\sigma (h'_1, h''_1)\sigma (h, h'_2h''_2)$ for all $h, h', h''\in H$).
If $A$ is a left $H$-comodule algebra with comodule structure
$a\mapsto a_{(-1)}\otimes a_{(0)}$, one may consider the new associative
product on $A$ given by $a*b=\sigma (a_{(-1)}, b_{(-1)})a_{(0)}b_{(0)}$.
This product is afforded by the map $T:A\otimes A\rightarrow A\otimes A$,
$T(a\otimes b)=\sigma (a_{(-1)}, b_{(-1)})a_{(0)}\otimes b_{(0)}$,
which is a pseudotwistor with companions $\tilde{T}_1$, $\tilde{T}_2$
given by the formulae
\begin{eqnarray*}
&&\tilde{T}_1(a\otimes b\otimes c)=
\sigma ^{-1}(a_{(-1)_1}, b_{(-1)_1})\sigma (a_{(-1)_2}, b_{(-1)_2}
c_{(-1)})a_{(0)}\otimes b_{(0)}\otimes c_{(0)}, \\
&&\tilde{T}_2(a\otimes b\otimes c)=\sigma ^{-1}(b_{(-1)_1}, c_{(-1)_1})
\sigma (a_{(-1)}b_{(-1)_2}, c_{(-1)_2})a_{(0)}\otimes b_{(0)}\otimes
c_{(0)}.
\end{eqnarray*}
In particular, for $A=H$, we obtain that the ``twisted bialgebra''
$_{\sigma }H$, with multiplication $a*b=\sigma (a_1, b_1)a_2b_2$,
for all $a, b\in H$, is obtained as a deformation of $H$ through the
pseudotwistor $T(a\otimes b)=\sigma (a_1, b_1)a_2\otimes b_2$ with companions
$\tilde{T}_1(a\otimes b\otimes c)=\sigma ^{-1}(a_1, b_1)\sigma (a_2,
b_2c_1)a_3\otimes b_3\otimes c_2$ and $\tilde{T}_2(a\otimes b\otimes c)=
\sigma ^{-1}(b_1, c_1)\sigma (a_1b_2, c_2)a_2\otimes b_3\otimes c_3$,
for all $a, b, c\in H$.}
\end{examples}

\begin{lemma} \label{com}
Let ${\mathcal{C}}$ be a (strict) braided monoidal category with
braiding $c$. Let $V$ be an object in ${\mathcal{C}}$ and
$T:V\otimes V \rightarrow V\otimes V$ a morphism in ${\mathcal{C}}$.
Then
\begin{eqnarray}
&&(V\otimes c_{V, V})\circ (T\otimes V)\circ (V\otimes c_{V, V}^{-1})=
(c_{V, V}^{-1}\otimes V)\circ (V\otimes T)\circ (c_{V, V}\otimes V),
\label{T1} \\
&&(V\otimes c_{V, V}^{-1})\circ (T\otimes V)\circ (V\otimes c_{V, V})=
(c_{V, V}\otimes V)\circ (V\otimes T)\circ (c_{V, V}^{-1}\otimes V),
\label{T2}
\end{eqnarray}
as morphisms $V\otimes V\otimes V\rightarrow V\otimes V\otimes V$ in
${\mathcal{C}}$. These two morphisms will be denoted by
$\tilde{T}_1(c)$ and $\tilde{T}_2(c)$ and will be called the
companions of $T$ with respect to the braiding $c$. If $c_{V,
V}^{-1}=c_{V, V}$ (for instance if ${\mathcal{C}}$ is symmetric),
the two companions coincide and will be simply denoted by
$T_{13}(c)$.
\end{lemma}

\begin{proof} The naturality of $c$ implies
$(V\otimes T)\circ c_{V\otimes V, V}=c_{V\otimes V, V}\circ (T\otimes V)$.
Since $c$ is a braiding, we have $c_{V\otimes V, V}=(c_{V, V}\otimes V)\circ
(V\otimes c_{V, V})$, hence we obtain
\begin{eqnarray*}
(V\otimes T)\circ (c_{V, V}\otimes V)\circ (V\otimes c_{V, V})=
(c_{V, V}\otimes V)\circ (V\otimes c_{V, V})\circ (T\otimes V).
\end{eqnarray*}
By composing to the left with $c_{V, V}^{-1}\otimes V$ and to the right with
$V\otimes c_{V, V}^{-1}$, we obtain the desired equality (\ref{T1}).
Similarly one can check that (\ref{T2}) holds, too.
\end{proof}

\begin{definition} \label{brtw}
Let ${\mathcal{C}}$ be a (strict) braided monoidal category, $(A,
\mu , u)$ an algebra in ${\mathcal{C}}$ and $T:A\otimes A\rightarrow
A\otimes A$ a morphism in ${\mathcal{C}}$. Assume that $c_{A,
A}^{-1}=c_{A, A}$ (so we have the morphism $T_{13}(c)$ in
${\mathcal{C}}$ as above). If $T$ is a pseudotwistor with companions
$\tilde{T}_1=\tilde{T}_2=T_{13}(c)$ and moreover $(T\otimes A)\circ
(A\otimes T)=(A\otimes T)\circ (T\otimes A)$, we will call $T$ a
\dtext{braided twistor} for $A$ in ${\mathcal{C}}$.
\end{definition}

Consider now ${\mathcal{C}}$ to be the category of
$\mathbb{Z}_2$-graded vector spaces, which is braided (even
symmetric) with braiding given by $c(v\otimes
w)=(-1)^{|v||w|}w\otimes v$, for $v, w$ homogeneous elements. If
$(\Omega , d)$ is a DG algebra, then $\Omega $ becomes a
$\mathbb{Z}_2$-graded algebra (i.e. an algebra in ${\mathcal{C}}$)
by putting even components in degree zero and odd components in
degree one. The map $T$ given by (\ref{fedtw}) is obviously a
morphism in ${\mathcal{C}}$, and using the above braiding one can
see that the morphism $T_{13}(c)$ in ${\mathcal{C}}$ is given by the
formula $T_{13}(c)(\omega \otimes \zeta \otimes \eta )=\omega
\otimes \zeta \otimes \eta -(-1)^{|\omega |+|\zeta |}d(\omega
)\otimes \zeta \otimes d(\eta )$, for homogeneous $\omega $, $\zeta
$, $\eta $ (which is {\it different} from the ordinary $T_{13}$),
and one can now check that $T$ is a braided twistor for $\Omega $ in
${\mathcal{C}}$, and obviously $\Omega ^T$ is just $\Omega $ endowed
with the Fedosov product, regarded as a $\mathbb{Z}_2$-graded
algebra.

\begin{theorem}\label{general}
Let $(A, \mu , u)$ be an algebra in a (strict) monoidal category
${\mathcal{C}}$, let $T, R:A\otimes A\rightarrow A\otimes A$ be
morphisms in ${\mathcal{C}}$, such that $R$ is an isomorphism and a
twisting map between $A$ and itself. Consider the morphisms
\begin{eqnarray}
&&\tilde{T}_1(R):=(R^{-1}\otimes A)\circ (A\otimes T)\circ (R\otimes A),
\label{T1R} \\
&&\tilde{T}_2(R):=(A\otimes R^{-1})\circ (T\otimes A)\circ (A\otimes R).
\label{T2R}
\end{eqnarray}
Define the morphism $P:=R\circ T:A\otimes A\rightarrow A\otimes A$.
Then:

\noindent\textbf{(i)} The relation (\ref{tw2}) holds for $P$ if and
only if (\ref{catw1}) holds for $T$, with $\tilde{T}_1(R)$ in place
of $\tilde{T}_1$.

\noindent\textbf{(ii)} The relation (\ref{tw3}) holds for $P$ if and
only if (\ref{catw2})
holds for $T$, with $\tilde{T}_2(R)$ in place of $\tilde{T}_2$. \\[2mm]
In particular, it follows that if $T$ is a pseudotwistor for $A$
with companions $\tilde{T}_1(R)$ and $\tilde{T}_2(R)$, then $P$ is a
twisting map between $A$ and itself.

\noindent\textbf{(iii)} Conversely, assume that $P$ is a twisting
map and the following relations are satisfied:
\begin{eqnarray}
&&(P\otimes A)\circ (A\otimes P)\circ (P\otimes A)=(A\otimes P)\circ
(P\otimes A)\circ (A\otimes P), \label{co1} \\
&&(R\otimes A)\circ (A\otimes R)\circ (R\otimes A)=(A\otimes R)\circ
(R\otimes A)\circ (A\otimes R), \label{co2} \\
&&(P\otimes A)\circ (A\otimes P)\circ (R\otimes A)=(A\otimes R)\circ
(P\otimes A)\circ (A\otimes P), \label{co3} \\
&&(R\otimes A)\circ (A\otimes P)\circ (P\otimes A)=(A\otimes P)\circ
(P\otimes A)\circ (A\otimes R). \label{co4}
\end{eqnarray}
Then $T$ is a pseudotwistor for $A$ with companions $\tilde{T}_1(R)$
and $\tilde{T}_2(R)$.

\noindent\textbf{(iv)} Assume that (iii) holds and moreover
\begin{eqnarray}
&&(P\otimes A)\circ (A\otimes R)\circ (R\otimes A)=(A\otimes R)\circ
(R\otimes A)\circ (A\otimes P), \label{co5} \\
&&(R\otimes A)\circ (A\otimes R)\circ (P\otimes A)=(A\otimes P)\circ
(R\otimes A)\circ (A\otimes R). \label{co6}
\end{eqnarray}
Then $R$ is also a twisting map between $A^T$ and itself.
\end{theorem}

\begin{proof} We prove (i), while (ii) is similar and left to the reader.
Assume first that (\ref{tw2}) holds for $P$. Then we can compute:
\begin{eqnarray*}
T\circ (A\otimes \mu )&=&R^{-1}\circ P\circ (A\otimes \mu )\\
&{{\rm (\ref{tw2})}\atop =}&R^{-1}\circ (\mu \otimes A)\circ (A\otimes P)
\circ (P\otimes A)\\
&=&R^{-1}\circ (\mu \otimes A)\circ (A\otimes R)\circ (A\otimes T)\circ
(R\otimes A)\circ (T\otimes A)\\
&{{\rm (\ref{tw2})}\atop =}&(A\otimes \mu )\circ (R^{-1}\otimes A)
\circ (A\otimes T)\circ (R\otimes A)\circ (T\otimes A)\\
&=&(A\otimes \mu )\circ \tilde{T}_1(R)\circ (T\otimes A),
\end{eqnarray*}
which is precisely the condition (\ref{catw1}). Conversely, assuming that
(\ref{catw1}) holds, we compute:
\begin{eqnarray*}
P\circ (A\otimes \mu )&=&R\circ T\circ (A\otimes \mu )\\
&{{\rm (\ref{catw1})}\atop =}&R\circ (A\otimes \mu )\circ \tilde{T}_1(R)\circ
(T\otimes A)\\
&=&R\circ (A\otimes \mu )\circ (R^{-1}\otimes A)\circ (A\otimes T)\circ
(R\otimes A)\circ (T\otimes A)\\
&{{\rm (\ref{tw2})}\atop =}&(\mu \otimes A)\circ (A\otimes R)
\circ (A\otimes T)\circ (R\otimes A)\circ (T\otimes A)\\
&=&(\mu \otimes A)\circ (A\otimes P)\circ (P\otimes A),
\end{eqnarray*}
which is (\ref{tw2}) for $P$.
Now we prove (iii). By (i) and (ii), it is enough to check (\ref{catw3}).
We compute:
\begin{eqnarray*}
\tilde{T}_1(R)\circ (T\otimes A)\circ (A\otimes T)&=&
(R^{-1}\otimes A)\circ (A\otimes T)\circ (R\otimes A)\circ (T\otimes A)
\circ (A\otimes T)\\
&=&(R^{-1}\otimes A)\circ (A\otimes R^{-1})\circ (A\otimes P)\circ
(P\otimes A)\circ \\
&&(A\otimes R^{-1})\circ (A\otimes P)\\
&{{\rm (\ref{co4})}\atop =}&(R^{-1}\otimes A)\circ (A\otimes R^{-1})
\circ (R^{-1}\otimes A)\circ (A\otimes P)\circ \\
&&(P\otimes A)\circ (A\otimes P)\\
&{{\rm (\ref{co1},\; \ref{co2})}\atop =}&(A\otimes R^{-1})
\circ (R^{-1}\otimes A)\circ (A\otimes R^{-1})\circ (P\otimes A)\circ \\
&&(A\otimes P)\circ (P\otimes A)\\
&{{\rm (\ref{co3})}\atop =}&(A\otimes R^{-1})\circ (R^{-1}\otimes A)\circ 
(P\otimes A)\circ (A\otimes P)\circ \\
&&(R^{-1}\otimes A)\circ (P\otimes A)\\
&=&(A\otimes R^{-1})\circ (T\otimes A)\circ (A\otimes R)\circ
(A\otimes T)\circ (T\otimes A)\\
&=&\tilde{T}_2(R)\circ (A\otimes T)\circ (T\otimes A).
\end{eqnarray*}
(iv) We check (\ref{tw2}) and leave (\ref{tw3}) to the reader. We compute:
\begin{eqnarray*}
R\circ (A\otimes \mu \circ T)&=&R\circ (A\otimes \mu \circ R^{-1}\circ P)\\
&=&R\circ (A\otimes \mu )\circ (A\otimes R^{-1})\circ (A\otimes P)\\
&{{\rm (\ref{tw2})}\atop =}&(\mu \otimes A)\circ (A\otimes R)
\circ (R\otimes A)\circ (A\otimes R^{-1})\circ (A\otimes P)\\
&{{\rm (\ref{co2})}\atop =}&(\mu \otimes A)\circ (R^{-1}\otimes A)\circ 
(A\otimes R)\circ (R\otimes A)\circ (A\otimes P)\\
&{{\rm (\ref{co5})}\atop =}&(\mu \otimes A)\circ (R^{-1}\otimes A)
\circ (P\otimes A)\circ (A\otimes R)\circ (R\otimes A)\\
&=&(\mu \circ R^{-1}\circ P\otimes A)\circ (A\otimes R)\circ (R\otimes A)\\
&=&(\mu \circ T\otimes A)\circ (A\otimes R)\circ (R\otimes A),
\end{eqnarray*}
finishing the proof.
\end{proof}

Our motivating example for Theorem \ref{general} was provided
by the theory of braided quantum groups, a concept introduced by M.
Durdevich in \cite{durde1} as a generalization of the usual braided
groups (=Hopf algebras in braided categories, in Majid's
terminology), which in turn contains as examples some important
algebras such as braided and ordinary Clifford algebras, see
\cite{durde2}. If $G=(A, \mu , \Delta , \varepsilon , S, \sigma )$
is a braided quantum group (so $\sigma $ is a bijective twisting map
between $A$ and itself) and $n\in \mathbb{Z}$, Durdevich defined
some operators $\sigma _n:A\otimes A\rightarrow A\otimes A$ and
proved that the maps $\mu _n:A\otimes A\rightarrow A$, $\mu _n=\mu
\circ \sigma _n^{-1}\circ \sigma $, give new associative algebra
structures on $A$ (with the same unit). This result may be regarded
as a consequence of Theorem \ref{general}. Indeed, for any $n$, the
maps $R:=\sigma _n$ and $P:=\sigma $ satisfy the hypotheses of the
theorem, hence the map $T:=R^{-1}\circ P=\sigma _n^{-1}\circ \sigma
$ is a pseudotwistor for $A$, giving rise to the associative
multiplication $\mu _n$.

More generally, if $A$ is an algebra, Durdevich introduced the
concept of {\it braid system} over $A$, as being a collection
$\mathcal{F}$ of bijective twisting maps between $A$ and itself,
satisfying the condition
\begin{eqnarray*}
&&(\alpha \otimes A)\circ (A\otimes \beta )\circ (\gamma \otimes A)=
(A\otimes \gamma )\circ (\beta \otimes A)\circ (A\otimes \alpha ),
\;\;\; \forall \;\alpha , \beta , \gamma \in \mathcal{F}.
\end{eqnarray*}
If we take $\alpha , \beta \in \mathcal{F}$ and define $T:A\otimes
A\rightarrow A\otimes A$, $T:=\alpha ^{-1}\circ \beta $, by Theorem
\ref{general} we obtain that $T$ is a pseudotwistor for $A$, giving
rise to a new associative multiplication on $A$.

We record the following two easy consequences of Theorem
\ref{general}.
\begin{corollary} \label{braid}
Let ${\mathcal{C}}$ be a (strict) braided monoidal category with
braiding $c$, $(A, \mu , u)$ an algebra in ${\mathcal{C}}$ and
$T:A\otimes A\rightarrow A\otimes A$ a morphism in ${\mathcal{C}}$;
assume also that $c_{A, A}^{-1}=c_{A, A}$. Define the morphism
$R:A\otimes A\rightarrow A\otimes A$ by $R:=c_{A, A}\circ T$. Then
$T$ satisfies the condition (\ref{catw1}) (respectively
(\ref{catw2})) with $T_{13}(c)$ in place of $\tilde{T}_1$
(respectively $\tilde{T}_2$) if and only if $R$ satisfies
(\ref{tw2}) (respectively (\ref{tw3})). In particular, if $T$ is a
braided twistor for $A$ in ${\mathcal{C}}$, then $R$ is a twisting
map between $A$ and itself.
\end{corollary}

\begin{corollary}
Let ${\mathcal{C}}$ be a (strict) braided monoidal category with
braiding $c$ and $(A, \mu , u)$ an algebra in ${\mathcal{C}}$. Then
$T:=c_{A, A}^2$ is a pseudotwistor for $A$ in ${\mathcal{C}}$ (this
follows by taking $R=c_{A, A}^{-1}$ and $P=c_{A, A}$ in Theorem
\ref{general}). In particular it follows that $(A, \mu \circ c_{A,
A}^2, u)$ is a new algebra in ${\mathcal{C}}$.
\end{corollary}

This algebra $(A, \mu \circ c_{A, A}^2, u)$ 
allows us to give an interpretation of the concept of
\dtext{ribbon algebra} introduced by Akrami and Majid in \cite{am},
as an essential ingredient for constructing braided Hochschild and
cyclic cohomology. Recall from \cite{am} that a ribbon algebra in a
braided monoidal category $({\mathcal{C}}, c)$ is an algebra $(A,
\mu , u)$ in ${\mathcal{C}}$ equipped with an isomorphism $\sigma
:A\rightarrow A$ in ${\mathcal{C}}$ such that $\mu \circ (\sigma
\otimes \sigma )\circ c_{A, A}^2= \sigma \circ \mu $ and $\sigma
\circ u=u$ (such a $\sigma $ is called a \dtext{ribbon automorphism}
for $A$). The naturality of $c$ implies $(\sigma \otimes \sigma
)\circ c_{A, A}^2=c_{A, A}^2\circ (\sigma \otimes \sigma )$, so the
above relation may be written as $\mu \circ c_{A, A}^2 \circ (\sigma
\otimes \sigma )=\sigma \circ \mu $. Hence, a ribbon automorphism
for $A$ is the same thing as an algebra isomorphism from $(A, \mu ,
u)$ to $(A, \mu \circ c_{A, A}^2, u)$.

Let $D$ be an algebra and $T$ a twistor for $D$. We intend to lift
$T$ to the algebra $\O D$ of universal differential forms on $D$; it
will turn out that the natural way of doing this does \textit{not}
provide a twistor, but a braided twistor. In order to simplify the
proof, we will use a braiding notation. Namely, we denote a braided
twistor $T$ for an algebra $A$ in a braided monoidal category with
braiding $c$ satisfying $c_{A,A}^{-1}=c_{A,A}$ by
\[
    \xy
        (0,0)*{\LTwistor{T}};
    \endxy
\]
where we will omit the label $T$ whenever there is no risk of confusion.
With this notation, the conditions for $T$ to be a braided twistor are
written as:

\[
    \includegraphics[scale=1]{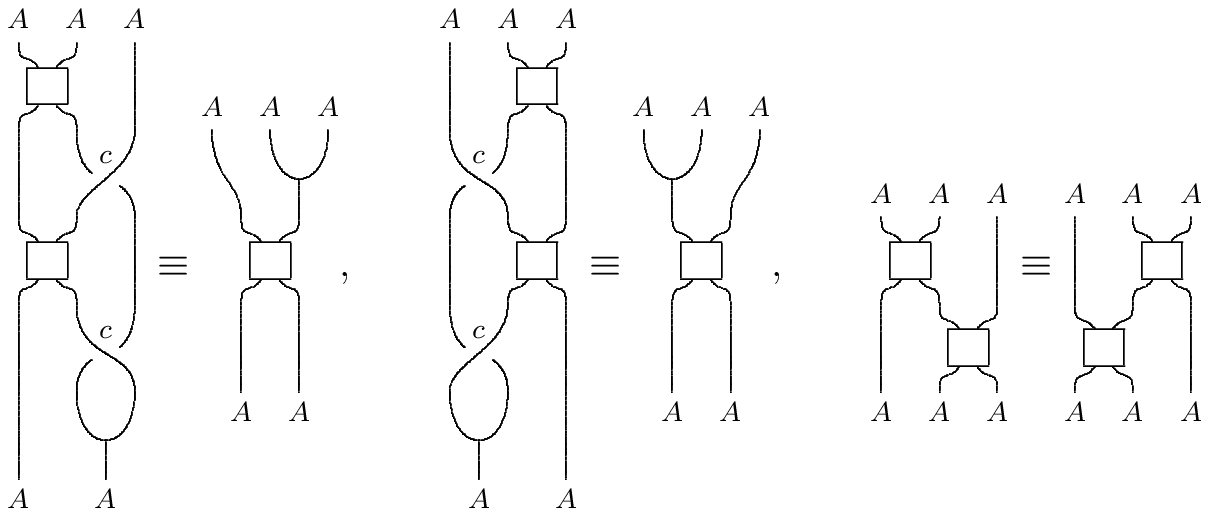}
\]

It is also worth writing the two equivalent definitions of $T_{13}(c)$
using this notation, namely:

\[
    \includegraphics[scale=1]{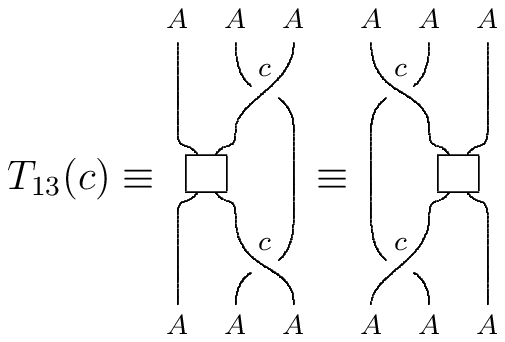}
\]

Let us consider now an algebra $D$ together with a twistor $T:D\ot
D\to D\ot D$. From Corollary \ref{braid} we know that the map
$R:=\tau\circ T$ is a twisting map between $D$ and itself. But then,
using Theorem \ref{twistdifferentialforms}, we may lift the twisting
map $R$ to a twisting map $\widetilde{R}:\O D\ot \O D\to \O D\ot \O
D$ between the algebra of universal differential forms $\O D$ and
itself. Using again Corollary \ref{braid} in the category of graded
vector spaces (with the graded flip $\tau_{gr}$ as a braiding) we
obtain that the map $\widetilde{T}:\O D\ot \O D\to \O D\ot \O D$
defined as $\widetilde{T}:=\tau_{gr}\circ \widetilde{R}$ satisfies
the conditions \eqref{catw1} and \eqref{catw2} with
$\tilde{T}_1=\tilde{T}_2=T_{13}(\tau_{gr})$. Moreover, it is clear
that $\widetilde{T}^0\equiv T$, since $\widetilde{R}$ extends $R$
and the graded flip coincides with the classical flip on degree 0
elements. Let us check that $\widetilde{T}$ also satisfies the
condition
\begin{equation}\label{catw5}
    (\widetilde{T}\ot \O D)\circ (\O D\ot \widetilde{T}) = 
(\O D\ot \widetilde{T})\circ (\widetilde{T}\ot \O D),
\end{equation}
and hence $\widetilde{T}$ is a braided (graded) twistor for the algebra $\O D$.
In order to do this, we follow a standard procedure when dealing 
with differential
calculi. First, as the restriction of $\widetilde{T}$ to $\O^0 D$ is a twistor,
it satisfies the condition. Second, assume that the condition is 
satisfied for an
element $\omega\ot \eta \ot \theta$ in $\O D\ot \O D\ot \O D$, 
and let us prove that
it is also satisfied for $d\omega\ot \eta \ot \theta$, 
$\omega\ot d\eta \ot \theta$
and $\omega\ot \eta \ot d\theta$. First of all, realize that, for homogeneous
$\omega,\eta\in \O D$, we have
\begin{equation}
    \tau_{gr}(\eta\ot d\omega) = (-1)^{|d\omega||\eta|}d\omega\ot \eta=
    (-1)^{(|\omega|+1) |\eta|}d\omega\ot \eta=
     (\varepsilon\ot d)\circ\tau_{gr}(\eta\ot \omega),
\end{equation}
where $d$ and $\varepsilon$ denote respectively the differential and
the grading of $\O D$. As a consequence of this equality and the
compatibilities of $\widetilde{R}$ with the differential (cf.
\eqref{twisteddiff1} and \eqref{twisteddiff2}), we realize
immediately that the map $\widetilde{T}$ satisfies the following
compatibility relations with the differential:
\begin{gather}
\widetilde{T}\circ (d\ot \O D) = (d\ot \O D)\circ \widetilde{T},
\label{twistordiff1}\\
\widetilde{T}\circ (\O D\ot d) = (\O D\ot d)\circ \widetilde{T}
\label{twistordiff2}.
\end{gather}

Using braiding knotation we have:

\[
    \includegraphics[scale=1]{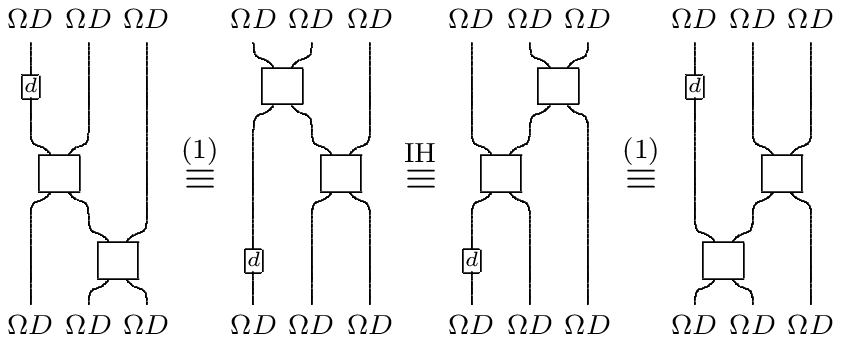}
\]
where in $(1)$ we are using \eqref{twistordiff1}, and in the second
equality we are using the induction hypothesis, and so the condition
\eqref{catw5} for $\widetilde{T}$ behaves well under the
differential in the first factor. The proof for the condition with
the differential on the second or third factors is similar, and left
to the reader.

Finally, we have to check that this condition also behaves well under 
products on any of the factors. For doing this, we need slightly 
stronger induction hypotheses. Namely, assume that we have 
$\omega_1,\omega_2,\eta,\theta$ such that the condition is satisfied for 
$\omega_i\ot \eta'\ot \theta'$, being $\eta',\theta'$ \textbf{any} elements 
in $\O D$ such that $|\eta'|\leq |\eta|$ and $|\theta'|\leq |\theta|$, 
i.e. we assume that the condition is true when we fix the $\omega_i$'s 
and let the $\eta'$ and $\theta'$ vary up to some degree bound, and 
let us prove that in this case the condition holds for 
$\omega_1\omega_2\ot \eta'\ot \theta'$. For this, take into account that 
$\widetilde{T}$ preserves the degree of homogeneous elements, since both 
$\widetilde{R}$ and $\tau_{gr}$ do. Now, we have

\[
    \includegraphics[scale=1]{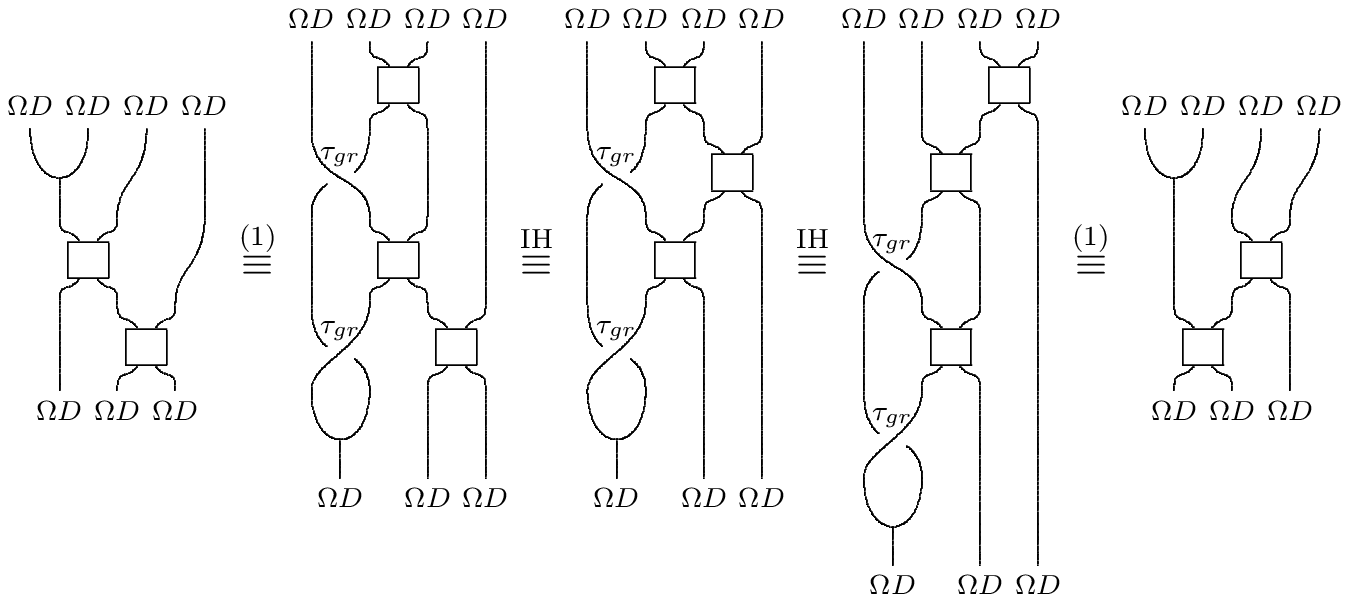}
\]
where in $(1)$ we are using \eqref{catw1}, in the equalities labeled
with IH we are using our strengthened induction hypotheses. The
desired result follows. Similar proofs exist when applying
multiplication in the second or third factors. It is easy to see
that, as a consequence of the properties we have just proved, we
obtain that the map $\widetilde{T}$ is a braided (graded) twistor on
the differential graded algebra $\O D$. More concretely, we have
proved the first part of the following result:

\begin{theorem}
    Let $D$ be an algebra and $T:D\ot D\to D\ot D$ a twistor for $D$.
    Consider $R:=\tau\circ T$, the twisting map associated to $T$.
    Let $\widetilde{R}$ be the extension of $R$ to $\O D$, then the
    map $\widetilde{T}:=\tau_{gr}\circ \widetilde{R}$ is a braided (graded)
    twistor for $\O D$. Moreover, the algebra $(\O D)^ {\widetilde{T}}$ is a
    differential graded algebra with differential $d$.
\end{theorem}

\begin{proof}
    The only part left to prove is that the map $d$ is still a differential
    for the deformed algebra $(\O D)^{\widetilde{T}}$, but this is an easy
    consequence of the fact that both the differential $d$ and the grading
    $\varepsilon$ commute with the twistor $\widetilde{T}$.
\end{proof}

The deformed algebra $(\O D)^{\widetilde{T}}$ has, as the
$0^{\text{th}}$ degree component, the algebra $D^T$, and, whenever
$T$ is bijective, it is generated (as a graded differential algebra)
by $D^ T$, henceforth $(\O D)^{\widetilde{T}}$ is a differential
calculus over $D^T$. Thus, as a consequence of the Universal
Property for the algebra of universal differential forms, we may
conclude that $(\O D)^{\widetilde{T}}$ is a quotient of the graded
differential algebra $\O (D^T)$.


\end{document}